\documentclass[11pt]{article}
\usepackage{amsmath,amssymb,amsfonts,amsthm,graphicx}
\usepackage{hyperref}
\newtheorem{theorem}{Theorem}[section]
\newtheorem{proposition}[theorem]{Proposition}

\newtheorem{lemma}[theorem]{Lemma}

\newtheorem{remark}[theorem]{Remark}
\newtheorem{definition}[theorem]{Definition}
\usepackage[titletoc,title]{appendix}
\numberwithin{equation}{section}
\title{ERGODIC PROPERTIES OF $\theta$-EXPANSIONS AND A GAUSS-KUZMIN-TYPE PROBLEM}
\author{
    Dan Lascu\footnote{e-mail: lascudan@gmail.com.}, Florin Nicolae\footnote{e-mail: florin.nicolae@anmb.ro.}\nonumber \\
    \emph{\small Mircea cel Batran Naval Academy, 1 Fulgerului, 900218 Constanta,
    Romania} 
    }
\begin{document}
\maketitle
\thispagestyle{empty}
\begin{abstract}
A generalization of the regular continued fractions was given by Chakraborty and Rao \cite{CR-2003}. 
For the transformation which generates this expansion and its invariant measure, the Perron-Frobenius operator is given and studied.
For this expansion, we apply the method of Rockett and Sz\"usz \cite{RS-1992} and obtained the solution of its Gauss-Kuzmin-type problem.
\end{abstract}

\section{Introduction}
Chakraborty and Rao \cite{CR-2003} considered a continued fraction expansion of a number in terms of an irrational $\theta \in (0,1)$. This new expansion of positive reals, different from the regular continued fraction expansion, is called $\theta$-expansion.

The purpose of this paper is to give some ergodic properties and to solve a Gauss-Kuzmin problem for $\theta$-expansions. 
In order to solve the Gauss-Kuzmin problem, we apply the method of Rockett and Sz\"usz \cite{RS-1992}.
First we outline the historical framework of this problem. 
In Section 1.2, we present the current framework. 
In Section 1.3, we review known results. 
In Section 1.4, the main theorem will be shown. 

\subsection{Gauss' Problem and its progress}

The development of the metric theory of continued fractions started in 1928 when Kuzmin proved a theorem due to Gauss. This theorem is now called Gauss-Kuzmin theorem.

Any irrational $0 \leq x < 1$ can be written as the infinite regular continued fraction 
\begin{equation}
x = \displaystyle \frac{1}{a_1+\displaystyle \frac{1}{a_2+\displaystyle \frac{1}{a_3+ \ddots}}} :=[0;a_1, a_2, a_3, \ldots], \label{1.1}
\end{equation}
where $a_n \in \mathbb{N}_+ : = \left\{1, 2, 3, \ldots\right\}$ \cite{IK-2002}.
Define the \textit{regular continued fraction} (or \textit{Gauss}) \textit{transformation} $\tau$ on the unit interval $I:=[0, 1]$ by 
\begin{equation}
\tau (x) = \left\{\begin{array}{lll} 
\displaystyle \frac{1}{x}-\left\lfloor \displaystyle \frac{1}{x} \right\rfloor & \hbox{if} & x \neq 0, \\ 
\\
0 & \hbox{if} & x = 0,
\end{array} \right. \label{1.2}
\end{equation}
where $\left\lfloor \cdot \right\rfloor$ denotes the floor (or entire) function.
With respect to the asymptotic behavior of iterations $\tau^{n}=\tau\circ \cdots\circ \tau$ ($n$-times) of $\tau$, in 1800 Gauss wrote (in modern notation) that
\begin{equation}
\lim_{n \rightarrow \infty} \lambda \left(\tau^n \leq x\right) = \frac{\log (1+x)}{\log 2},
\quad x \in I, \label{1.3}
\end{equation}
where $\lambda$ denotes the Lebesgue measure on $I$.
In 1812, Gauss asked Laplace \cite{Brez} to estimate the {\it $n$-th error term} $e_n(x)$ defined by
\begin{equation}
e_n(x) := \lambda (\tau^{-n}[0, x]) - \frac{\log (1+x)}{\log 2}, \quad n \geq 1, \ x\in I. \label{1.4}
\end{equation}
This has been called \textit{Gauss' Problem}. 
In 1928 Kuzmin \cite{Kuzmin-1928} showed that $e_n(x) = \mathcal{O}(q^{\sqrt{n}})$ as $n \rightarrow \infty$, uniformly in $x$ with some (unspecified) $0 < q < 1$. 
Independently, L\'evy \cite{Levy-1929} proved in 1929 that $\left|e_n(x)\right| \leq q^n$ for $n \in \mathbb{N}_+$, $x \in I$, with $q = 0.67157...$. 
For such historical reasons, the {\it Gauss-Kuzmin-L\'evy theorem} 
is regarded as the first basic result in the rich metrical theory of continued fractions. 

Apart from the regular continued fraction expansion, very many other continued fraction expansions were studied \cite{RS-1992,Schweiger}. By such a development, generalizations of these problems for non-regular continued fractions 
are also called as the {\it Gauss-Kuzmin problem} and the {\it Gauss-Kuzmin-L\'evy problem} \cite{IS-2006, L-2013, L-2014, Sebe-2000, Sebe-2001, Sebe-2002, SL-2014}.

\subsection{$\theta$-expansions as dynamical system}

In this paper, we consider a generalization of the Gauss transformation and prove an analogous result. 

This transformation was studied by Chakraborty, Rao and Dasgupta in \cite{CR-2003} and \cite{CD-2004} and by Sebe and Lascu in \cite{SL-2014}.

Fix an irrational $ \theta \in (0,1)$. In \cite{CR-2003}, Chakraborty and Rao showed that any $x \in \left(0, \theta \right)$ can be written as the form
\begin{equation}
x = \frac{1}{\displaystyle a_1\theta
+\frac{1}{\displaystyle a_2\theta 
+ \frac{1}{\displaystyle a_3\theta + \ddots} }} := [0; a_1 \theta, a_2 \theta, \ldots], \label{1.5}
\end{equation}
where $a_n$'s are non-negative integers. 
We will simply write (\ref{1.5})
\begin{equation}
x := [a_1 \theta, a_2 \theta, \ldots]. \label{1.6}
\end{equation}

Such $a_n$'s are called {\it incomplete quotients (or continued fraction digits)} 
of $x$ with respect to the expansion in (\ref{1.5}) in this paper.

This continued fraction is treated as the following dynamical systems.

\begin{definition} \label{def.1.1}
Let $\theta \in (0,1)$ and $m \in \mathbb{N}_+$ such that $\theta^2 = 1/m$. 
\begin{enumerate}

\item [(i)]
The measure-theoretical dynamical system $([0,\theta],{\cal B}_{[0,\theta]}, T_{\theta})$ is defined as follows:
$\mathcal{B}_{[0,\theta]}$ denotes the $\sigma$-algebra of all Borel subsets of $[0,\theta]$, and $T_{\theta}$ is the transformation 
\begin{equation}
T_{\theta}: [0,\theta] \to [0,\theta];\quad 
T_{\theta}(x):=
\left\{
\begin{array}{ll}
{\displaystyle \frac{1}{x} - \theta \left \lfloor \frac{1}{x \theta} \right\rfloor} & 
{\displaystyle \hbox{if } x \in (0, \theta],}\\
\\
0 & \hbox{if } x=0.
\end{array}
\right. \label{1.7}
\end{equation} 

\item [(ii)] 
In addition to (i), we write $([0,\theta], {\cal B}_{[0,\theta]}, \gamma_{\theta}, T_{\theta} )$ as $( [0,\theta],{\cal B}_{[0,\theta]},T_{\theta})$ with the following probability measure $\gamma_{\theta}$ on $( [0,\theta], {\cal B}_{[0,\theta]})$:

\begin{equation}
\gamma_{\theta} (A) := 
\frac{1}{\log \left(1+\theta^{2}\right)}
\int_{A} \frac{\theta dx}{1 + \theta x}, 
\quad A \in {\mathcal{B}}_{[0,\theta]}. \label{1.8}
\end{equation}

\end{enumerate}
\end{definition}

Define the {\it quantized index map} $\eta : [0, \theta] \to {\mathbb N}$ by
\begin{equation}
\eta(x) := \left\{\begin{array}{ll} 
\left\lfloor \displaystyle \frac{1}{x \theta}\right\rfloor  & \hbox{if }  x \neq 0, \\ 
\\
\infty & \hbox{if }  x = 0.
\end{array} \right. \label{1.9}
\end{equation}
By using $T_{\theta}$ and $\eta$, the sequence $(a_{n})_{n \in \mathbb N_+}$ in (\ref{1.5}) is obtained as follows:
\begin{equation}
a_n(x) = \eta(T_{\theta}^{n-1}(x)), \quad n \geq 1, \label{1.10}
\end{equation}
with $T_{\theta}^0 (x) = x$.

In this way, $T_{\theta}$ gives the algorithm of $\theta$-expansion.

\begin{proposition} \label{prop.1.2}
Let 
$([0,\theta],{\cal B}_{[0,\theta]},\gamma_{\theta},T_{\theta})$ be as in Definition \ref{def.1.1}(ii).
\begin{enumerate}
\item [(i)]
$([0,\theta],{\cal B}_{[0,\theta]},\gamma_{\theta},T_{\theta})$ is ergodic.
\item [(ii)]
The measure $\gamma_{\theta}$ is invariant under $T_{\theta}$, that is, 
$\gamma_{\theta} (A) = \gamma_{\theta} (T^{-1}_{\theta}(A))$ 
for any $A \in {\mathcal{B}}_{[0, \theta]}$.
\end{enumerate}
\end{proposition}
\noindent \textbf{Proof.} 
See Section 8 in \cite{CR-2003}.
\hfill $\Box$\\

By Proposition \ref{prop.1.2}(ii), $([0,\theta],{\cal B}_{[0,\theta]},\gamma_{\theta},T_{\theta})$ is a 
``dynamical system" in the sense of Definition 3.1.3 in \cite{BG-1997}.

\subsection{Known results and applications}

For $\theta$-expansions, we show known results and their applications in this subsection.

Let $0 < \theta < 1$ and $m \in \mathbb{N}_+$ such that $\theta^2 = 1/m$. 
In what follows the stated identities hold for all $n$ in case $x$ has an infinite $\theta$-expansion and they hold for $n \leq k$ in case $x$ has a finite $\theta$-expansion terminating at the $k$-th stage \cite{CR-2003}.

To this end, define real functions $p_n(x)$ and $q_n(x)$, for $n \in \mathbb{N}_+$, by 
\begin{eqnarray}
p_n(x) &:=& a_n(x) \theta p_{n-1}(x) + p_{n-2}(x), \quad \label{1.11} \\
q_n(x) &:=& a_n(x) \theta q_{n-1}(x) + q_{n-2}(x), \quad \label{1.12}
\end{eqnarray}
with $p_{-1}(x) := 1$, $p_0(x) := 0$, $q_{-1}(x) := 0$ and $q_{0}(x) := 1$. 
By using (\ref{1.11}) and (\ref{1.12}), we can verify that 
\begin{equation}
x = \frac{p_n(x) + T^n_{\theta}(x)p_{n-1}(x)}
{q_n(x) + T^n_{\theta}(x)q_{n-1}(x)}, \quad n \geq 1 \label{1.13}
\end{equation}
and 
\begin{equation}
x - \frac{p_n(x)}{q_n(x)} = \frac{(-1)^{n+1} T^n_{\theta}(x)}{q_n(x)(q_n(x)+T^n_{\theta}(x)q_{n-1}(x))}, \quad n \geq 1. \label{1.14}
\end{equation}

We introduced a partition of the interval $[0, \theta]$ which is natural to the $\theta$-expansions. Such a partition is generated by the \textit{fundamental intervals} (or \textit{cylinders}) of rank $n$. 
For any $n \in \mathbb{N}_+$ and $i^{(n)}=(i_1, \ldots, i_n) \in \mathbb{N}_m^n$, define the {\it fundamental interval associated with} $i^{(n)}$ by
\begin{equation}
I \left(i^{(n)}\right) = \{x \in [0, \theta]:  a_k(x) = i_k \mbox{ for } k=1, \ldots, n \}, \label{2.01}
\end{equation}
where $I \left(i^{(0)}\right) = [0, \theta]$. 
We will write $I(a_1, \ldots, a_n) = I\left(a^{(n)}\right)$, $n \in \mathbb{N}_+$. 
If $n \geq 1$ and $i_n \in \mathbb{N}_m$, then we have
\begin{equation}
I(a_1, \ldots, a_n) = I\left(i^{(n)}\right). \label{2.03}
\end{equation}

Using the ergodicity of $T_{\theta}$ and Birkhoff's ergodic theorem \cite{DK-2002}, a number of results were obtained.

For $q_n$ in (\ref{1.11}), its asymptotic growth rate $\beta$ is defined as 
\begin{equation}
\beta = \lim_{n \rightarrow \infty} \frac{1}{n} \log q_n. \label{1.15}
\end{equation}
This is a L\'evy result.  
Although the calculation algorithm is correct, Chakraborty and Rao \cite{CR-2003} misspelled the expression of $\beta$. 
They should write
\begin{equation}
\beta = \frac{-1}{\log(1 + \theta^2)} \int_{0}^{\theta} \frac{\theta \log x}{1 + \theta x}\mathrm{d}x. \label{1.16}
\end{equation}
They also give a Khintchin result, i.e., the asymptotic value of the arithmetic mean of $a_1, a_2, \ldots, a_n$ where $a_n$'s are given in (\ref{1.10}). We have
\begin{equation}
\lim_{n \rightarrow \infty} \frac{a_1 + a_2 + \ldots + a_n}{n} = \infty. \label{1.17}
\end{equation}

In \cite{SL-2014}, Sebe and Lascu proved a Gauss-Kuzmin theorem for the transformation $T_{\theta}$. In order to solve the problem, they applied the theory of random systems with complete connections (RSCC) by Iosifescu \cite{IG-2009}.
We remind that a random system with complete connections is a quadruple
\begin{equation}
\left\{\left([0, \theta], {\mathcal{B}}_{[0, \theta]}\right), \left(\mathbb{N}_m, {\mathcal{P}}(\mathbb{N}_m)\right), u, P\right\}, \label{1.18}
\end{equation}
where $u : [0, \theta] \rightarrow [0, \theta]$, 
\begin{equation}
u(x) = u_i(x) := \frac{1}{i \theta + x}, \label{1.19}
\end{equation}
and $P$ is a transition probability function from $\left([0,\theta],{\mathcal{B}}_{[0,\theta]}\right)$ to $\left(\mathbb{N}_m, {\mathcal{P}}(\mathbb{N}_m) \right)$ given by 
\begin{equation}
P(x) = P_i(x) := \frac{\theta x + 1}{(x + i\theta)(x + (i+1)\theta )}. \label{1.20}
\end{equation}
Here $\mathbb{N}_m := \{m, m+1, \ldots\}$,
and $\mathcal{P} \left(\mathbb{N}_m\right)$ denotes the power set of $\mathbb{N}_m$.
Also, the associated Markov operator of RSCC (\ref{1.18}) is denoted by $U$ and has the transition probability function 
\begin{equation}
Q(x, A) = \sum_{i \in W(x,A)}P_{i}(x), \quad x \in [0,\theta], \ A \in {\mathcal{B}}_{[0, \theta]}, \label{1.21}
\end{equation}
where $W(x,A) = \left\{i \in \mathbb{N}_m: u_{i}(x) \in A \right\}$.

Using the asymptotic and ergodic properties of operators associated with RSCC (\ref{1.18}), i.e., the ergodicity of RSCC, they obtained a convergence rate result for the Gauss-Kuzmin-type problem.

For more details about using RSCC in solving Gauss-Kuzmin-L\'evy-type theorems, see \cite{L-2013,Sebe-2000,Sebe-2001,Sebe-2002,SL-2014}

\subsection{Main theorem}

We show our main theorem in this subsection. 

The measure $\gamma_{\theta}$ in (\ref{1.8}) is the unique absolutely continuous invariant measure for the map $T_{\theta}$ in (\ref{1.7}). 
In particular, if one iterates any other absolutely continuous invariant measure repeatedly by $T_{\theta}$, 
it will converge exponentially to $\gamma_{\theta}$.

Let $\mu$ be a non-atomic probability measure on $\mathcal{B}_{[0,\theta]}$ and define 
\begin{eqnarray}
F_n (x) &:=& \mu (T_{\theta}^n \leq x), \ x \in [0,\theta], \ n \in \mathbb{N}_+, \label{1.29} \\
F(x) &:=& \displaystyle \lim_{n \rightarrow \infty}F_n(x), \ x \in I, \label{1.30}
\end{eqnarray}
with $F_0 (x) = \mu ([0,x])$. 

Then the following holds.
\begin{theorem} [A Gauss-Kuzmin-type theorem] \label{G-K}
Let $T_{\theta}$, $\gamma_{\theta}$ and $F_n$ be as in $\mathrm{(\ref{1.7})}$, $\mathrm{(\ref{1.8})}$ and $\mathrm{(\ref{1.29})}$, respectively. 
Then there exists a constant $0 < q < \theta$ such that $F_n$ is written as 
\begin{equation}
F_n (x) = \frac{\log ((m\theta + x)\theta)}{\log (1+\theta^2)} + \mathcal{O}(q^n). \label{1.31}
\end{equation}
\end{theorem} 
\begin{remark}
\rm{
From (\ref{1.31}), we see that
\begin{equation}
F(x) = \gamma_{\theta}([0,x]). \label{1.32}
\end{equation}
In fact, the Gauss-Kuzmin theorem estimates the error 
\begin{equation}
e_{\theta}(x) := e_{\theta}(x, \mu) = \mu (T_{\theta}^n \leq x) - \gamma_{\theta} ([0,x]), \quad x \in [0,\theta]. \label{1.33}
\end{equation}
}
\end{remark}
\noindent \textbf{Open problem.}
Solve the Gauss-Kuzmin-L\'evy problem of $T_{\theta}$.
For example, study the optimality of the convergence rate by using the same strategy as in \cite{IS-2006}.

The rest of the paper is organized as follows.
In Section 2, we give some ergodic properties of $\theta$-expansions 
and we find the entropy of the transformation which generates $\theta$-expansion. 
In Section 3, we derive the associated Perron-Frobenius operator under different probability measures on $([0, \theta], {\mathcal{B}}_{[0, \theta]})$. 
We treat the Perron-Frobenius operator of $([0, \theta],{\cal B}_{[0, \theta]},\gamma_{\theta}, T_{\theta})$, 
and derive its asymptotic behavior.
In Section 4, we prove Theorem \ref{G-K} for $\theta$-expansions. 
In Section 4.1, we give the necessary results used to prove the Gauss-Kuzmin theorem. 
The essential argument of the proof is the Gauss-Kuzmin-type equation. 
We will also give some results concerning the behavior of the derivative of $\{F_{n}\}$ in (\ref{1.29}) which will allow us to complete the proof of the Theorem \ref{G-K} in Section 4.2.

\section{Ergodic properties and entropy} 

Let $\theta \in (0,1)$ and $m \in \mathbb{N}_+$ such that $\theta^2 = 1/m$. 
In this section we give some consequences of ergodicity in terms of properties of the continued fraction expansion for almost every real number $x \in (0, \theta)$. 

The following theorem presents another L\'evy and Khintchin-type results. 
The relation (\ref{1.15}) allows to find other asymptotical results of this type. 
Applying our machinery to the ergodic system $([0, \theta], {\mathcal{B}}_{[0, \theta]}, \gamma_{\theta}, T_{\theta})$ 
we obtain a result on the statistical data of $\theta$-expansions, i.e., 
almost sure asymptotics for the geometrical mean value for partial quotients. 
\begin{proposition} \label{prop.1.3}

For allmost all $x \in (0, \theta)$ one has
\begin{equation}
\lim_{n \to \infty} \frac{1}{n} \log \left(\lambda_{\theta} \left(I\left(i^{(n)}\right)\right)\right) = -2 \beta, \qquad \qquad \qquad \ \ \label{1.22}
\end{equation}
\begin{equation}
\lim_{n \to \infty} \frac{1}{n} \log \left|x - \frac{p_n}{q_n}\right| = -2 \beta, \qquad \qquad \qquad \qquad \quad \ \  \label{1.23} 
\end{equation}
\begin{equation}
\lim_{n \to \infty} (a_1 a_2 \cdots a_n)^{1/n} = \prod^{\infty}_{k=m} \left(1+\frac{1}{k(k+2)}\right)^{\frac{\log k}{\log(1+\theta^2)}}. \label{1.24} 
\end{equation}
Here, $\lambda_{\theta}$ denote a Lebesgue measure on $[0, \theta]$ and $I\left(i^{(n)}\right)$ are the cylinders in $\mathrm{(\ref{2.01})}$, and $a_n$'s, $p_n$ and $q_n$ are given in $\mathrm{(\ref{1.10})}$, $\mathrm{(\ref{1.11})}$ and $\mathrm{(\ref{1.12})}$, respectively.
\end{proposition}
\noindent \textbf{Proof.} See Appendix. 
\hfill $\Box$

As it is well known, entropy is an important concept of information in physics, chemistry, and information theory \cite{PY-1998}. It can be seen as a measure for the amount of ``disorder" of a system. Entropy also plays an important role in ergodic theory. Like Birkhoff's ergodic theorem \cite{DK-2002,PY-1998} the entropy is a fundamental result in ergodic theory.
For a measure preserving transformation, its entropy is often defined by using partitions. We choose this partition with respect to the fundamental intervals $I\left(i^{(n)}\right)$ in (\ref{2.01}). The entropy $h(T_{\theta})$ is obtained from L\'evy's result (\ref{1.22}) and the theorem of Shannon-McMillan-Breiman \cite{DK-2002}:
\begin{equation}
h(T_{\theta}) = \lim_{n \to \infty} -\frac{1}{n} \log \left(\lambda_{\theta} \left(I\left(i^{(n)}\right)\right)\right) = 2 \beta.
\end{equation}
One could also compute entropy by the beautiful formula of Rohlin \cite{Roh-1964}:
\begin{eqnarray}
h(T_{\theta}) &=& \int^{\theta}_{0} \log\left|T'_{\theta} (x)\right| d \gamma_{\theta}(x) \nonumber \\
              &=& \int^{\theta}_{0} \frac{- \log x^2}{\log(1+\theta^2)}\frac{\theta dx}{1+\theta x} = 2 \beta, \label{1.28}
\end{eqnarray}
where $\gamma_{\theta}$ and $\beta$ are given in (\ref{1.8}) and (\ref{1.16}), respectively.

\section{The Perron-Frobenius operator of $T_{\theta}$ under $\gamma_{\theta}$}

Let $([0, \theta],{\cal B}_{[0, \theta]},\gamma_{\theta}, T_{\theta})$ be as in Definition 1.1. 
In this section, we derive its Perron-Frobenius operator.

Let $\mu$ be a probability measure on $([0, \theta], {\mathcal{B}}_{[0, \theta]})$ 
such that $\mu((T_{\theta})^{-1}(A)) = 0$ whenever $\mu(A) = 0$ for $A \in {\mathcal{B}}_{[0, \theta]}$. 
For example, this condition is satisfied if $T_{\theta}$ is $\mu$-preserving, that is, $\mu (T_{\theta})^{-1} = \mu$. 
Let %
$
L^1([0, \theta], \mu):=\{f: [0, \theta] \rightarrow \mathbb{C} : \int^{\theta}_0 |f |\mathrm{d}\mu < \infty \}.
$
The {\it Perron-Frobenius operator} $U$ of $([0, \theta],{\cal B}_{[0, \theta]},\mu, T_{\theta})$ 
is defined as the bounded linear operator on the Banach space $L^1([0, \theta],\mu)$ 
such that the following holds:
\begin{equation}
\int_{A}Uf \,\mathrm{d}\mu = \int_{(T_{\theta})^{-1}(A)}f\, \mathrm{d}\mu \quad 
\mbox{ for all }
A \in {\mathcal{B}}_{[0, \theta]},\, f \in L^1([0, \theta],\mu). \label{5.29}
\end{equation}
About more details, see \cite{BG-1997, IK-2002} or Appendix A in \cite{L-2013}.
\begin{proposition} \label{prop-PF}
Let $([0, \theta],{\cal B}_{[0, \theta]},\gamma_{\theta}, T_{\theta})$ be as in Definition 1.1, and let $U$ denote its Perron-Frobenius operator. 
\begin{enumerate}
\item[(i)] 
The following equation holds: 
\begin{equation}
Uf(x) = \sum_{i \geq m}P_{i}(x)\,f(u_{i}(x)), \quad m \in \mathbb{N}_+, \, 
f \in L^1([0, \theta],\gamma_{\theta}), \label{5.30}
\end{equation}
where $P_{i}$ and $u_i$, $i \geq m$, are as in $\mathrm{(\ref{1.20})}$ and $\mathrm{(\ref{1.19})}$, respectively.
\item[(ii)]
Let $\mu$ be a probability measure on ${\mathcal{B}}_{[0,\theta]}$. 
Assume that $\mu$ is absolutely continuous with respect to the Lebesgue measure $\lambda_{\theta}$ (and denote $\mu \ll \lambda_{\theta}$, i.e., if $\mu(A) = 0$ for every set $A$ with $\lambda_{\theta}(A) = 0$) and let $h = d\mu / d \lambda_{\theta}$  \mbox{a.e. in } $[0, \theta]$. 
Then the following holds:
\begin{enumerate}
	\item[(a)] The Perron-Frobenius operator $S$ of $T_{\theta}$ under $\mu$ is given a.e. in $[0, \theta]$ by the equation 
	\begin{eqnarray}
S f(x) &=& \frac{1}{h(x)} \sum_{i \geq m} \frac{1}{i \theta + x} f(u_i(x)) \label{5.31} \\
             &=& \frac{U g(x)}{(1 + \theta x)h(x)}, \ f \in L^1([0, \theta], \mu), \label{5.32}
\end{eqnarray}
where $g(x):=(1 + \theta x)f(x)h(x)$, $x \in [0, \theta]$.
In addition, the $n$-th power $S^n$ of $S$ is given as follows:
\begin{equation}
S^n f(x) = \frac{U^n g(x)}{(1 + \theta x)h(x)} \label{5.33}
\end{equation}
for any $f \in L^1([0, \theta], \mu)$ and any $n \in \mathbb{N}_+$.
  \item[(b)]
The Perron-Frobenius operator $V$ of \, $T_{\theta}$ under $\lambda_{\theta}$ is given a.e. in $[0, \theta]$ by the equation 
\begin{equation}
V f(x) = \sum_{i \geq m} \frac{1}{(i \theta + x)^2} f(u_i(x)), \ f \in L^1([0, \theta], \lambda_{\theta}). \label{5.34}
\end{equation}
The powers of $V$ are given a.e. in $[0, \theta]$ by the equation
\begin{equation}
V^n f(x) = \frac{U^n g(x)}{1 + \theta x}, \ f \in L^1([0, \theta], \lambda_{\theta}), \mbox{ } n \in \mathbb{N}_+, \label{5.35}
\end{equation}
where $g(x):=(1 + \theta x)f(x)$, $x \in [0, \theta]$. 
  \item[(c)]
For any $n \in \mathbb{N}_+$ and $A \in {\mathcal{B}}_{[0,\theta]}$,
we have 
\begin{equation}
\mu \left((T_{\theta})^{-n}(A)\right) 
= \int_{A} U^nf(x) \mathrm{d} \gamma_{\theta}(x), \label{2.11}
\end{equation}
where $f(x):= (\log(1+\theta^2)) \frac{1 + x \theta}{\theta ^2}h(x)$, $x \in [0, \theta]$. 
\end{enumerate}
\end{enumerate}
\end{proposition}

\noindent \textbf{Proof.} 
See Appendix. 
\hfill $\Box$
\\

For a function $f : [0,\theta] \to {\Bbb C}$, define the {\it variation} ${\rm var}_{A}f$ of $f$ on a subset $A$ of $[0, \theta]$ by
\begin{equation}
{\rm var}_A f := \sup \sum^{k-1}_{i=1} |f(t_i) - f(t_{i-1})|, \label{5.36}
\end{equation}
where the supremum being taken over $t_1 < \cdots < t_k$, $t_i \in A$, 
$1 \leq i \leq k$, and $k \geq 2$. 
We write simply $\mathrm{var} f$ for $\mathrm{var}_{[0, \theta]} f$.
Let $L^{\infty}([0,\theta])$ denote the collection of all bounded measurable functions 
$f: [0,\theta] \rightarrow \mathbb{C}$. 
It is known that $L^{\infty}([0,\theta]) \subset L^1([0,\theta])$. 
Let  $L([0,\theta])$ denote the Banach space of all complex-valued Lipschitz continuous functions on $[0,\theta]$
with the following norm $\|\cdot\|_L$:
\begin{equation}
\left\| f \right\|_L := \sup_{x \in [0,\theta]} |f(x)| + s(f), \label{5.36'}
\end{equation}
with 
\begin{equation}
s(f):=\sup_{x \ne y} \frac{|f(x) - f(y)|}{|x - y|}, \quad f \in L([0,\theta]). \label{5.36''}
\end{equation}

In the following proposition we show that the operator $U$ in (\ref{5.30}) preserves monotonicity and 
enjoys a contraction property for Lipschitz continuous functions.

\begin{proposition} \label{prop-var} 
Let $U$ be as in $(\ref{5.30})$.  
\begin{enumerate}
  \item[(i)]
  Let $f \in L^{\infty}([0,\theta])$. Then the following holds:
\begin{enumerate}
    \item[(a)]
	If $f$ is non-decreasing (non-increasing), 
	then $Uf$ is non-increasing (non-decreasing).
    \item[(b)]
	If $f$ is monotone, then
\begin{equation}
\mathrm{var}\,(Uf) \leq k_m \cdot \mathrm{var } f \quad \mbox{where } k_m := \frac{1}{m+1}. \label{5.37}
\end{equation}
\end{enumerate}

 \item[(ii)] 
 For any $f \in L([0, \theta])$, we have
\begin{equation}
s(Uf) \leq q \cdot s(f), \label{5.38}
\end{equation}
where
\begin{equation}
q:=m \left( \sum_{i \geq m} \left(\frac{m}{i^3(i+1)}+ \frac{i+1-m}{i(i+1)^3}\right) \right) \label{5.38'}
\end{equation}
\end{enumerate}
\end{proposition} 

\noindent \textbf{Proof.} See Appendix. 
\hfill $\Box$

\section{Proof of Theorem \ref{G-K}} 

In this section we prove our main theorem applying the method of Rockett and Sz\"usz \cite{RS-1992}.
Let $\theta \in (0,1)$ and $m \in \mathbb{N}_+$ such that $\theta^2 = 1/m$.

\subsection{Necessary lemmas} 

In this subsection, we show some lemmas. First, we show that $\{F_n\}$ in (\ref{1.29}) satisfy a Gauss-Kuzmin-type equation.

\begin{lemma} \label{G-K-eq}
For functions $\{F_n\}$ in $\mathrm{(\ref{1.29})}$, the following Gauss-Kuzmin-type equation holds:
\begin{equation}
F_{n+1} (x) = \sum_{i \geq m }\left\{F_{n}\left(\frac{1}{i \theta}\right) - F_{n}\left(\frac{1}{i \theta + x}\right)\right\} \label{2.1}
\end{equation}
for $x \in [0,\theta]$ and $n \in \mathbb{N}$.
\end{lemma}

\noindent\textbf{Proof.}
Let $I_{n}=\{x \in [0,\theta]: T^n_{\theta}(x) \leq x \}$ and 
\begin{equation}
I_{n,i} = \left\{ x \in I_{n}: \ \frac{1}{i \theta + x} <  T^n_{\theta}(x) < \frac{1}{i \theta} \right\}. \label{2.2}
\end{equation}
From (\ref{1.7}) and (\ref{1.10}), we see that
\begin{equation}
T_{\theta}^{n}(x) = \displaystyle\frac{1}{a_{n+1} \theta + T_{\theta}^{n+1}(x)}, \quad n \in \mathbb{N}_+. \label{2.3}
\end{equation}
From the definition of $I_{n,i}$ and (\ref{2.3}) it follows that for any $n \in \mathbb{N}$, $I_{n+1} = \bigcup_{i \geq m } I_{n,i}$. 
From this and using the $\sigma$-additivity of $\mu$, we have
\begin{equation}
\mu(I_{m,n+1})=\mu\left(\displaystyle \bigcup_{i \in \mathbb{N}} I_{m,n,i}\right)=\sum_{i \in \mathbb{N}}\mu(I_{m,n,i}). \label{2.4}
\end{equation}
Then (\ref{2.1}) holds because $F_{n+1} (x)= \mu(I_{n+1})$ and 
\begin{equation}
\mu(I_{n,i}) = F_{n}\left(\frac{1}{i \theta}\right) - F_{n}\left(\frac{1}{i \theta + x}\right). \label{2.5}
\end{equation}
\hfill $\Box$
\begin{remark}
\rm{
Assume that for some $p \in \mathbb{N}$, the derivative $F'_p$ exists everywhere in $[0, \theta]$ and is bounded. 
Then it is easy to see by induction that $F'_{p+n}$ exists and is bounded for all $n \in \mathbb{N}_+$. 
This allows us to differentiate (\ref{2.1}) term by term, obtaining 
\begin{equation}
F'_{n+1}(x) = \sum_{i \geq m} \frac{1}{(i \theta + x)^2}F'_{n}\left(\frac{1}{i \theta + x}\right). \label{2.6}
\end{equation}
}
\end{remark}

\noindent We introduce functions $\{ f_{n} \}$ as follows:
\begin{equation}
f_{n}(x) := (1+ \theta x)F'_{n}(x), \quad x \in [0, \theta], \ n \in \mathbb{N}. \label{2.7}
\end{equation}
Then (\ref{2.6}) is 
\begin{equation}
f_{n+1}(x) = \sum_{i \geq m} P_i(x)f_{n}\left(u_i(x)\right), \label{2.8}
\end{equation}
where $P_i(x)$ and $u_i(x)$ are given in (\ref{1.20}) and (\ref{1.19}), respectively. By Proposition \ref{prop-PF} (i), we have that $f_{n+1}(x) = Uf_n(x)$.
\begin{lemma} \label{lema-2.3}
For $\{f_{n}\}$ in $\mathrm{(\ref{2.7})}$, define $M_n := \displaystyle \max_{x \in [0, \theta]}|f'_{n}(x)|$. 
Then
\begin{equation}
M_{n+1} \leq q \cdot M_n, \label{2.9}
\end{equation}
where $q$ is the constant in $\mathrm{(\ref{5.38'})}$.

\end{lemma}
\noindent\textbf{Proof.}
Since 
\begin{equation}
P_i(x) = \frac{1}{\theta} \left[ \frac{1-i \theta^2}{x + i \theta} - \frac{1 - (i+1) \theta^2}{x + i \theta} \right], \label{2.11}
\end{equation}
we have
%
%
%
\begin{equation}
f'_{n+1}(x) = \sum_{i \geq m} \frac{1-(i+1)\theta^2}{(x+i\theta)(x+(i+1)\theta)^3} f'_n(\alpha_i) - \sum_{i \geq m} \frac{P_i(x)}{(x+i\theta)^2} f'_n(u_i(x)), \label{2.12}
\end{equation}
where $u_{i+1}(x) < \alpha_i < u_i(x)$. Now (\ref{2.12}) implies 
\begin{equation}
M_{n+1} \leq M_{n} \cdot \max_{x \in [0, \theta]} \left( \sum_{i \geq m} \frac{(i+1)\theta^2 - 1}{(x+i \theta)(x + (i+1)\theta)^3} 
+ \sum_{i \geq m} \frac{P_i(x)}{(x+i\theta)^2} \right). \label{2.13}
\end{equation}
We now must calculate the maximum value of the sums in this expression. 
First, we note that
\begin{equation}
\frac{(i+1)\theta^2 - 1}{(x+i \theta)(x + (i+1)\theta)^3} \leq \frac{(i+1)\theta^2 - 1}{i(i+1)^3}, \label{2.14}
\end{equation}
where we use that $0 \leq x \leq \theta$.
Next, let 
\begin{equation}
h(x) := \sum_{i \geq m} \frac{P_i(x)}{(x+i\theta)^2}. \label{2.15}
\end{equation}
By Proposition \ref{prop-PF} (i) and Proposition \ref{prop-var} (i)(a), we have that the function $h$ is decreasing for $x \in [0, \theta]$ and $i \geq m$. Hence, $h(x) \leq h(0)$. This leads to
\begin{equation}
\sum_{i \geq m} \frac{P_i(x)}{(x+i\theta)^2} \leq m^2 \cdot \sum_{i \geq m} \frac{1}{i^3(i+1)}. \label{2.16}
\end{equation}
The relations (\ref{2.13}), (\ref{2.14})  and (\ref{2.16}) imply (\ref{2.9}) and that $q$ is as in (\ref{5.38'}). 

$\hfill \Box$

\subsection{Proof of Theorem \ref{G-K}}

Introduce a function $R_{n}(x)$ such that
\begin{equation}
F_{n}(x) = \frac{\log ((m\theta + x)\theta)}{\log (1+\theta^2)} + R_{n}(x). \label{2.17}
\end{equation}

Because $F_{n}(0)=0$ and $F_{n}(\theta)=1$, we have $R_{n}(0)=R_{n}(\theta)=0$. To prove Theorem \ref{G-K}, we have to show the existence of a constant $0 < q < \theta$ such that 
\begin{equation}
R_{n}(x) = {\mathcal O}(q^n). \label{2.18}
\end{equation}

If we can show that $f_{n}(x) = \frac{1}{\log (1+\theta^2)} + {\mathcal O}(q^n)$, then its integration will show the equation (\ref{1.31}).

To demonstrate that $f_{n}(x)$ has this desired form, it suffices to prove the following lemma.
\begin{lemma} \label{lema-2.4}
For any $x \in [0, \theta]$ and $n \in \mathbb{N}$, there exists a constant $0 < q < \theta$ such that 
\begin{equation}
f'_{n}(x) = {\mathcal O}(q^n). \label{2.19}
\end{equation}
\end{lemma}

\noindent \textbf{Proof}.
Let $q$ be as in (\ref{5.38'}). Using Lemma \ref{lema-2.4}, to show (\ref{2.19}) it is enough to prove that $q < \theta$ which is an easy task if we use a mathematical software. 
For example, for $m := 10$ we have $\theta=0.316228$ and $q=0.0533201$, and for $m := 17$ we have $\theta=0.242536$ and $q=0.0305636$.
$\hfill \Box$
%



\begin{appendices}

\section{Proofs of propositions}

We prove Propositions \ref{prop.1.3}, \ref{prop-PF} and \ref{prop-var} in this section.
\\

\subsection{Proof of Proposition \ref{prop.1.3}}

(i)
In \cite{SL-2014} it was shown that 
\begin{equation}
\lambda_{\theta}\left(I\left(a^{(n)}\right)\right) = \frac{1}{q_n (q_n + \theta q_{n-1})}. \label{A1}
\end{equation}

Thus
\begin{equation}
\frac{1}{(1+\theta)q_{n}^{2}} < \lambda_{\theta}\left(I\left(a^{(n)}\right)\right) < \frac{1}{q_{n}^{2}}, \label{A2}
\end{equation}
or
\begin{equation}
-\log(1+\theta) - 2 \log(q_n) < \lambda_{\theta}\left(I\left(a^{(n)}\right)\right) < -2 \log(q_n). \label{A3}
\end{equation}
Now apply (\ref{1.15}) to obtain (\ref{1.22}). 
\hfill $\Box$

(ii) In \cite{SL-2014} it was shown that
\begin{equation}
\frac{1}{q_n(q_{n+1}+\theta q_n)} \leq \left|x - \frac{p_n}{q_n}\right| \leq \frac{1}{q_nq_{n+1}}. \label{A4}
\end{equation}
Now (\ref{1.23}) follows from (\ref{1.15}) and {\ref{A4}}. 
\hfill $\Box$

(iii) Let $f(x):= \log(a_1(x))$ for $x \in (0, \theta)$, where $a_1$ is as in (\ref{1.10}). 
Then $f \in L^1\left((0, \theta), \gamma_{\theta}\right)$, i.e., $f$ is an integrable function, since
\begin{eqnarray}
\int_{0}^{\theta} f \mathrm{d}\gamma_{\theta} &=& 
\sum_{k=m}^{\infty} \int_{\frac{1}{(k+1)\theta}}^{\frac{1}{k\theta}} f \mathrm{d}\gamma_{\theta} = 
\sum_{k=m}^{\infty} \frac{1}{\log{(1+\theta^2)}} \int_{\frac{1}{(k+1)\theta}}^{\frac{1}{k\theta}} \frac{\theta \log a_1(x)}{1+ \theta x} \mathrm{d}x \nonumber \\
&=& \frac{1}{\log(1+\theta^2)} \sum_{k=m}^{\infty} \log k \cdot \log\left(1+\frac{1}{k(k+2)}\right) \nonumber \\
&=& \frac{1}{\log(1+\theta^2)} \sum_{k=m}^{\infty} \frac{\log k }{k(k+2)}. \label{A5}
\end{eqnarray}
Since the series $\sum_{k=m}^{\infty} (\log k)/(k(k+2))$ is convergent it follows that 
$\int_{0}^{\theta}f \mathrm{d}\gamma_{\theta} := s \in \mathbb{R}$.
Now we have
\begin{equation}
\lim_{n \to \infty} (a_1 a_2 \cdots a_n)^{1/n} = e^s = \prod^{\infty}_{k=m} \left(1+\frac{1}{k(k+2)}\right)^{\frac{\log k}{\log(1+\theta^2)}}. 
\end{equation}
\hfill $\Box$

\subsection{Proof of Proposition \ref{prop-PF}}

(i) See \cite{SL-2014}.

\noindent (ii)(a) 
Let $T_{\theta,i}$ denote the restriction of $T_{\theta}$ to the subinterval 
$I_i:=\left( \frac{1}{\theta (i+1)}, \frac{1}{\theta i}\right]$, $i \geq m$, $m \in \mathbb{N}$, that is, 
\begin{equation}
T_{\theta,i}(x) = \frac{1}{x} - \theta i, \quad x \in I_i. \label{A7}
\end{equation}
Let $C(A):=T^{-1}(A)$ and 
$C_{i}(A):=\left(T_{\theta,i}\right)^{-1}(A)$ for $A\in {\cal B}_{[0, \theta]}$.
Since $C(A)=\bigcup_{i}C_i(A)$ and $C_i\cap C_j$ is a null set when $i \ne j$,
we have
\begin{equation}
\int_{C(A)} f \,\mathrm{d} \gamma_{\theta} 
= \sum_{i \geq m} \int_{C_i(A)}f\, \mathrm{d} \gamma_{\theta},
\quad
f \in L^1([0, \theta],\gamma_{\theta}),\,A \in {\mathcal{B}}_{[0, \theta]}. \label{A8}
\end{equation}
From (\ref{A7}), for any $f \in L^1([0, \theta],\gamma_{\theta})$ 
and  $A \in {\mathcal{B}}_{[0, \theta]}$, we have
\begin{eqnarray}
\displaystyle \int_{C(A)} f(x) \,\mathrm{d}\mu(x) 
&=& \displaystyle \sum_{i \geq m} \displaystyle \int_{C_i(A)}f(x) \,\mathrm{d}\mu(x) \nonumber\\
\nonumber \\
& =& \displaystyle \sum_{i \geq m} \frac{1}{\theta} \displaystyle \int_{C_i(A)}f(x)h(x)\, \mathrm{d}x \nonumber \\
\nonumber \\
&=& \displaystyle \sum_{i \geq m} \frac{1}{\theta} \displaystyle \int_{A} \displaystyle \frac{f(u_{i}(y))\,h(u_{i}(y))}{(\theta i + y)^2} \,\mathrm{d}y \nonumber \\
\nonumber \\
&=& \frac{1}{\theta} \int_{A} \sum_{i \geq m} \displaystyle \frac{f(u_{i}(x))\,h(u_{i}(x))}{(\theta i + x)^2} \,\mathrm{d}x. \label{A9} 
\end{eqnarray}
Since $\mathrm{d}\mu = h \mathrm{d} \lambda_{\theta}$, (\ref{5.31}) follows from (\ref{A9}). 

Now, since $g(x)=(\theta x +1)f(x)h(x)$, from (\ref{5.30}) we have 
\begin{equation}
U g(x) = (\theta x +1) \sum_{i \geq m} \frac{h(u_{i}(x))}{(\theta i +x)^2}\,f(u_{i}(x)). \label{A10}
\end{equation}
Now, (\ref{5.32}) follows immediately from (\ref{5.31}) and (\ref{A10}). 
Using mathematical induction (\ref{5.33}) follows easily. 
\hfill $\Box$

(b) The formula (\ref{5.34}) is a consequence of (\ref{5.32}) and follows immediately. 

(c) See \cite{SL-2014}.
\hfill $\Box$

\subsection{Proof of Proposition \ref{prop-var}}

(i)(a)
To make a choice assume that $f$ is non-decreasing. 
Let $x < y$, \, $x, y \in [0, \theta]$. 
We have $U f(y) - U f(x) = S_1 + S_2$, where
\begin{eqnarray}
S_1 &=& \sum_{i \geq m} P_i(y) \left(f(u_i(y)) - f(u_i(x))\right), \label{A11} \\
S_2 &=& \sum_{i \geq m} \left(P_i(y) - P_i(x)\right)f(u_i(x)). \label{A12}
\end{eqnarray}
Clearly, $S_1 \leq 0$. Now, since $\sum_{i \geq m} P_i(x) = 1$ for any $x \in [0, \theta]$,
we can write 
\begin{equation}
S_2 = - \sum_{i \geq m} \left( f(u_m(x)) - f(u_i(x)) \right) \left(P_i(y) - P_i(x)\right). \label{A13}
\end{equation}
As is easy to see, the functions $P_i$ are increasing for all $i \geq m$.
Also, using that $f(u_m(x)) \geq f(u_i(x))$, we have that $S_2 \leq 0$. Thus $U f(y) - U f(x) \leq 0$ and the proof is complete.
\hfill $\Box$

(i)(b)
Assume that $f$ is non-decreasing. Then by (a) we have
\begin{equation}
\mathrm{var}\, Uf = Uf(0) - Uf(\theta) = \sum_{i \geq m} \left( P_i(0) f(u_i(0)) - P_i(\theta) f(u_i(\theta)) \right). \label{A14}
\end{equation}
By calculus, we have
\begin{eqnarray}
\mathrm{var}\, Uf &=& \sum_{i \geq m} \left( \frac{m}{i(i+1)}f\left(\frac{1}{\theta i}\right) - \frac{m+1}{(i+1)(i+2)}f\left(\frac{1}{\theta (i+1)}\right)\right) \nonumber\\
&=& \frac{1}{m+1} f(\theta) - \sum_{i \geq m} \frac{1}{(i+1)(i+2)}f\left(\frac{1}{\theta (i+1)}\right) \nonumber \\
&\leq& \frac{1}{m+1} f(\theta) - \sum_{i \geq m} \frac{1}{(i+1)(i+2)}f(0)  \nonumber \\
&=& \frac{1}{m+1} (f(\theta) - f(0)) = \frac{1}{m+1}\mathrm{var} f. \nonumber
\end{eqnarray}
\hfill $\Box$

(ii) For $x \ne y$, $x,y \in [0, \theta]$, we have
\begin{eqnarray}
\frac{Uf(y) - Uf(x)}{y - x} &=& \sum_{i \geq m} \frac{P_i(y) - P_i(x)}{y - x}f(u_i(x)) \nonumber \\
&-& \sum_{i \geq m} P_i(y) \frac{f(u_i(y)) - f(u_i(x))}{u_i(y) - u_i(x)}\cdot u_i(x)u_i(y). \label{A15}
\end{eqnarray}
Remark that 
\begin{equation}
P_i(x) = \frac{\theta (i+1-m)}{\theta (i+1) + x} + \frac{\theta (m-i)}{\theta i + x}, \quad i \geq m, \label{A16}
\end{equation}
and then
\begin{eqnarray}
&{}& \sum_{i \geq m} \frac{P_i(y) - P_i(x)}{y - x}f(u_i(x)) \nonumber \\
&=& \sum_{i \geq m} \frac{\theta (i+1-m)}{(y + \theta (i+1))(x + \theta (i+1))} \left( f(u_{i+1}(x)) - f(u_{i}(x)) \right). \label{A17}
\end{eqnarray}

Assume that $x > y$. It then follows from (\ref{A15}) and (\ref{A17}) that 
\begin{equation}
\left|\frac{Uf(y) - Uf(x)}{y - x}\right| \leq s(f) \sum_{i \geq m} \left( \frac{\theta^2(i+1-m)}{(y+\theta i)(y + \theta (i+1))^3} + \frac{P_i(y)}{(y+\theta i)^2}\right). \label{A18}
\end{equation}
Since the sum of right side of (\ref{A18}) is $q$ (see (\ref{2.13})) and since 

%
%
%
%
%
\begin{equation}
s(Uf) = \sup_{x \ne y} \left|\frac{Uf(y) - Uf(x)}{y - x}\right|
\end{equation}
then the proof is complete. 
\hfill $\Box$

\end{appendices}

\end{document}